\documentclass[plain]{aiaa-pretty}
\usepackage{amsthm}
\usepackage{amsmath}
\usepackage{amsfonts}
\usepackage{amssymb}
\usepackage{dsfont}
\usepackage{subfigure}
\usepackage{psfrag}
\usepackage{balance}
\usepackage{float}
\usepackage[titletoc,title]{appendix}
\usepackage{booktabs, makecell}
\usepackage{threeparttable}

\let\doi\relax
\usepackage{doi}

\newcommand{\overbar}[1]{\mkern 1.5mu\overline{\mkern-1.5mu#1\mkern-1.5mu}\mkern 1.5mu}
\newcommand{\virg}[1]{``#1''}

\newtheorem{remark}{Remark}

\author{Mirko Leomanni\thanks{Research Associate, Department of Information Engineering and Mathematics, \texttt{$\!$leomanni@diism.unisi.it}},
Gianni Bianchini\thanks{Professor, Department of Information Engineering and Mathematics, \texttt{giannibi@diism.unisi.it}},
Andrea Garulli \thanks{Professor, Department of Information Engineering and Mathematics, \texttt{garulli@diism.unisi.it}},
Renato Quartullo \thanks{PhD Student, Department of Information Engineering and Mathematics, \texttt{quartullo@diism.unisi.it}}\\[2mm]
\textit{Universit\`a di Siena, Siena, 53100, Italy}}

\title{Optimal Low-Thrust Orbit Transfers Made Easy:\\ A Direct Approach}

\abstract{
The optimization of low-thrust, multi-revolution orbit transfer trajectories is often regarded as a difficult problem in modern astrodynamics. In this paper, a flexible and computationally efficient approach is presented for the optimization of low-thrust orbit transfers under eclipse constraints. The proposed approach leverages a new dynamic model of the orbital motion and a Lyapunov-based initial guess generation scheme that is very easy to tune. A multi-objective, single-phase formulation of the optimal control problem is devised, which provides a convenient way to trade off fuel consumption and time of flight. A distinctive feature of such a formulation is that it requires no prior information about the structure of the optimal solution. Simulation results for two benchmark orbit transfer scenarios indicate that minimum-time, minimum-fuel and mixed time/fuel-optimal instances of the control problem can be readily solved via direct collocation, while incurring a significantly lower computational demand with respect to existing techniques.}

\begin{document}
\maketitle

\section*{Nomenclature}
\noindent\begin{tabular}{@{}lcl@{}}
$a,\,e,\,i,\,\omega,\,\Omega,\,v$&=&Classical orbital elements\\
$p,\,f,\,g,\,h,\,k,\,L$ &=&Modified equinoctial elements\\
$\rho,\,e_x,e_y,h_x,h_y,\sigma$&=&Ideal elements\\
$\mathbf{e}$&=&Relative eccentricity vector\\
$g_0$&=&Standard gravity, m/s$^2$\\
$I_{sp}$&=&Specific impulse, s\\
$m$&=&Spacecraft mass, kg\\
$\mathbf{q}$&=&Thrust unit vector\\
$R_e$&=&Earth radius, m\\
$R_s$&=&Sun radius, m\\
$\mathbf{r}_{es}$&=&Satellite-Earth vector, m\\
$\mathbf{r}_s$&=&Earth-Sun vector, m\\
$\mathbf{r}_{ss}$&=&Satellite-Sun vector, m\\
$s$&=&Ideal anomaly, rad\\
$T_{max}$&=&Maximum deliverable thrust, N\\
$t$&=&Time, s\\
$\mathbf{u}$&=&Control acceleration, m/s$^2$\\
$\mathbf{x}$&=&System state\\
$\overline{\mathbf{z}}$&=&Target equinoctial elements\\
$\alpha$&=&Time/fuel trade-off parameter\\
$\gamma$&=&Relative inclination, rad\\
$\epsilon$&=&Mesh error tolerance\\
$\Delta m$&=&Mass variation, kg\\
$\Delta t$&=&Time of flight, days\\
$\lambda_1,\,\lambda_2$&=&Nodal elements, rad\\
$\mu$&=&Earth gravitational parameter, m$^3$/s$^2$

\end{tabular}

\noindent\begin{tabular}{@{}lcl@{}}
$\xi$&=&Efficiency factor\\
$\xi_{cut}\hspace{2.2cm}$&=&Cut-off parameter\\
$\varrho_R$&=&Radial perturbation, m/s$^2$\\
$\varrho_T$&=&Transverse perturbation, m/s$^2$\\
$\varrho_N$&=&Normal perturbation, m/s$^2$\\
$\boldsymbol{\varrho}$&=&Perturbing acceleration vector, m/s$^2$\\
$\boldsymbol{\varrho}_{J2}$&=&J2 acceleration vector, m/s$^2$\\
$\zeta$&=&Throttle command\\
$\eta$&=&Throttle control input\\
$\theta_e$&=&Angular size of Earth radius, rad\\
$\theta_{es}$&=&Earth-Sun angle, rad\\
$\theta_s$&=&Angular size of Sun radius, rad\\
$\varphi_{az},\,\varphi_{el}$&=&Tuning parameters, rad\\
$\psi$&=&Shadow function\\
$\psi_\ell$&=&Smoothed shadow function\\[2mm]
\end{tabular}

\section{Introduction}
In recent years, low-thrust propulsion technologies such as Electric Propulsion (EP) have become popular as the primary means of propulsion for both planetary and interplanetary space missions. Their main appeal is a fuel efficiency which is ten times higher than that of conventional chemical thrusters. This feature is especially relevant in applications requiring large delta-v increments, a notable case being the orbital transfer about a central body. The low-thrust orbit transfer problem can be stated as the determination of a continuous spacecraft trajectory that satisfies initial and terminal conditions defined along two different orbits, while minimizing fuel expenditure and/or time of flight (TOF). Except for very special cases, this is a nonlinear optimal control problem (OCP) which is extremely hard to solve. Complications arise mainly due to the large time scale of the maneuver and to the bang-off-bang structure of the optimal control policy. An evident indication of such difficulties is that most papers on the subject deal only with a specific objective (either time or fuel optimization) or a specific transfer scenario.

Techniques for addressing the challenging issues in low-thrust trajectory optimization have been investigated for many years. Historically, the solution strategies have been divided
into two general categories: indirect methods and direct methods \cite{Bettssurvey98}. In indirect methods, necessary conditions for optimality are derived from optimal control theory via variational arguments. The necessary conditions form a boundary-value problem which is solved numerically in order to determine the optimal solution, see, e.g., \cite{Kechichian97,Haberkorn04,Topputo15}. The primary advantage of indirect methods is that they provide highly accurate solutions. Their main drawback is the very high sensitivity of the numerical solution to the initial guess for the costate variables, which makes these methods only suitable for systems of relatively low dimension. Direct methods, on the other hand, transcribe the continuous-time optimal control problem into a static nonlinear programming problem (NLP) which is solved by a nonlinear optimizer, see, e.g., \cite{Hargraves87,Enright92,betts2000}. This alleviates the sensitivity issue, as the solution procedure is much more robust with respect to the selection of the initial trajectory guess. In particular, direct pseudospectral methods (see, e.g., \cite{Elnagar95,Ross04,garg11}) have been applied to low-thrust problems with considerable success. In these methods, the state and control variables are parameterized using polynomials, and the system evolution is approximated via orthogonal collocation. Regardless of the transcription strategy, however, the application of direct methods typically requires the solution of large and computationally demanding NLPs. For a complete overview on these topics, the reader is referred to the dedicated books \cite{Betts10,conway2010,kechichian2018}.

Despite the advances in optimization techniques and tools, the low-thrust orbit transfer problem is still a computationally difficult one. Given the practical need to ease the computation, a number of studies have focused on approximate or suboptimal solution strategies. In \cite{Kluever98,Gao07}, control parametrization and averaging are exploited in order to reduce the size of the optimization problem. Although these approaches are more computationally efficient compared to those based on the exact dynamic model, the accumulated errors may be large, so that the true spacecraft trajectory may fail to reach the target orbit. Lyapunov-based guidance laws such as the Q-law \cite{petropoulos2004low,petropoulos2005refinements} have also been proposed in order to rapidly explore the solution space. These are usually quite robust to modeling errors, but lack optimality due to their heuristic nature. For these reasons, they are commonly employed as an initial guess for subsequent optimization (see, e.g., \cite{whiffen2006mystic}). In many cases, however, finding an adequate set of tuning parameters for such guidance laws can be as time consuming as the optimization process itself. Along a different line, some efforts have been made to study the impact of different parameterizations of the orbital motion on the computational efficiency of the optimization process, see, e.g., \cite{aziz2018low,Junkins19}. In particular, it is found in \cite{aziz2018low} that regularizing the orbital dynamics provides an effective way to speed-up the computation. In this paper, we will show that large computational gains are possible by working on the parametrization, while retaining the ability to solve the control problem to full optimality.

Another major source of complexity in low-thrust trajectory optimization is due to eclipsing. This is particularly relevant since most of the existing EP system are solar powered and therefore they are unable to operate during solar eclipse periods. Indeed, as pointed out in \cite{Betts15low,aziz2019smoothed}, low-thrust trajectory optimization without concerns for eclipsing is futile when subsequent analysis reveal that a thrusting maneuver has been scheduled in the shadow region. The net effect of solar eclipses is to introduce a state-dependent discontinuity in the control input structure, which constitutes a serious challenge for gradient-based optimization techniques. Unfortunately, this important issue is often disregarded in the literature. In the seminal paper \cite{ferrier2001}, the low-thrust orbit transfer problem with eclipsing is tackled by adopting an indirect method. More recently, the problem has been framed in the context of direct optimization. In particular, two different formulations have demonstrated their utility in applications: single-phase and multi-phase ones. The multi-phase formulation divides the OCP into multiple phases in which the system dynamics are smooth, while ensuring the continuity of the state trajectory among different phases through an appropriate set of event constraints. In this way, discontinuities are removed from the dynamic model and the optimization is cast over the switch times. A limitation of this method is that the structure of the optimal solution must be guessed a priori, which typically involves a rather complex machinery. For instance, it is not trivial to determine whether a phase should be added to the formulation.  Multi-phase approaches have been employed in \cite{Betts15low} and \cite{Rao16low} for the solution of minimum-fuel and minimum-time problems with eclipsing, respectively. The single-phase formulation is defined in the usual way and considers the OCP as a whole. An advantage of this method over multi-phase formulations is that it requires no prior information about the structure of the optimal solution, although some kind of smoothing is needed for the system dynamics. In \cite{aziz2019smoothed}, a single-phase approach is employed in combination with a smoothed eclipse model in order to solve a minimum-fuel problem involving a large number of orbital revolutions. This is done by using hybrid differential dynamic programming \cite{lantoine2012hybrid}. In \cite{Shannon20lowthrust}, the minimum-time problem with eclipsing is similarly cast and it is solved via collocation. A suboptimal solution to the minimum-fuel problem, which consists of applying the Q-law in the eclipsed part of the transfer and the optimal control policy in the remaining part, is also discussed.

The aim of this paper is to present a new approach for the optimization of low-thrust orbit transfers, employing a direct collocation method.
The work is in the same spirit of \cite{aziz2019smoothed,Shannon20lowthrust}, but some key features are introduced in the problem definition. First, equinoctial elements are replaced by a novel nonsingular set of \emph{ideal elements} inspired by Hansen's theory (see \cite{hansen1859,musen1960modified,deprit1975ideal}) for the description of the orbital motion. The rationale behind this choice is the observation that the dynamics of the \emph{ideal anomaly} do not directly depend on the perturbing acceleration, as opposed to those of all other anomalies. This makes the ideal anomaly attractive as a basis for regularization (see, e.g., \cite{Junkins04nonlinear, roa2016regularization}). Moreover, ideal elements enjoy a higher level of sparsity of the dynamic model structure compared to other element sets. This is advantageous because the OCP to be solved turns out to be a sparse NLP. As a second contribution, a flexible optimization framework is established, which unifies different types of performance requirements under the same single-phase OCP formulation. The adopted cost function involves a convex combination of TOF and fuel consumption, providing a suitable way to trade off these objectives by means of a scalar parameter. Eclipse effects are taken into account by adapting the smoothing technique in \cite{aziz2019smoothed} to the the proposed state parametrization. The OCP is transcribed into a sparse NLP by using the commercial package GPOPS–II \cite{patterson2014gpops}, which implements an \emph{hp}-adaptive Legendre–Gauss–Radau pseudospectral collocation strategy. The NLP is then solved by using a sparse nonlinear optimizer. Furthermore, a new Lyapunov-based guidance scheme that addresses the initial guess generation problem is devised. The guidance scheme features only three assignable parameters, which can be tuned in an intuitive manner.

The proposed optimization architecture is tested on two orbit transfer scenarios taken from the literature: a GTO-GEO transfer and a LEO-GEO transfer. Simulation results show that the method is able to solve minimum-time, mimimum-fuel and mixed time/fuel-optimal problem instances involving many eclipse transitions in a computationally efficient manner. Moreover, the capability of the method to solve large-scale optimization problems on low-power hardware is demonstrated for a realization featuring nearly half a million NLP variables. Overall, the obtained results indicate that the adopted parametrization and optimization scheme allow one to tackle the low-thrust OCP with relative ease in comparison to previous approaches.

The paper is organized as follows. Section~\ref{sec2} illustrates the new parametrization of the orbital motion and Section~\ref{sysmod} details the dynamic model used for trajectory optimization as well as the eclipse smoothing technique. Section~\ref{sec4} discusses the proposed OCP formulation and Section~\ref{initGG} presents the Lyapunov guidance scheme employed for the initial guess generation. Detailed simulation case studies of the optimization architecture are presented in Section~\ref{numsim} for the considered orbit transfer scenarios. Section~\ref{sec7} summarizes the main findings of this research and finalizes the paper.

\section{Parametrization of the Orbital Motion}\label{sec2}
The first step towards the derivation of an optimal orbit transfer strategy is to define a suitable set of parameters describing the orbital motion and to determine how such parameters evolve in response to perturbations. In this section, the perturbed Kepler problem is briefly reviewed and a new parametrization of the orbital motion is presented for this problem. The proposed parametrization is based on a nonsingular set of orbital parameters inspired by Hansen's theory.
\subsection{Perturbed Kepler Problem}
The perturbed Kepler problem amounts to describing the evolution of an orbit in response to a perturbing acceleration that accounts for all contributions other than point mass gravity. The classical solution to this problem is expressed in terms of the orbital elements $\{a,\,e,\,i,\,\omega,\,\Omega,\,v\}$, namely the semimajor axis, eccentricity, inclination, argument of periapsis, right ascension of the ascending node, and true anomaly. Let $\varrho_R$, $\varrho_T$ and $\varrho_N$ be the radial, transverse and normal components of the perturbing acceleration, expressed in the Radial-Transverse-Normal frame centered at the satellite. Whenever nonconservative perturbations such as thrusting are involved, it is customary to describe the evolution of the classical elements through Gauss' variational equations
\begin{equation}\label{gvecoe}
\setlength\arraycolsep{2pt}
\begin{array}{lll}
\dot{a}&=&\dfrac{2 a^2}{\sqrt{\mu a(1-e^2)}}\Big[e\sin (v) \,\varrho_R+(1+e\cos (v))\,\varrho_T\Big]\vspace{2mm} \\
\dot{e}&=&\sqrt{\dfrac{a(1-e^2)}{\mu}}\left[\sin (v) \,\varrho_R+\dfrac{e+(2+e\cos (v))\cos (v)}{(1+e\cos (v))}\varrho_T\right]\vspace{2mm} \\
\dot{i}&=&\sqrt{\dfrac{a(1-e^2)}{\mu}}\left[\dfrac{\cos(v+\omega)}{1+e\cos (v)}\varrho_{N}\right]\vspace{2mm} \\
\dot{\Omega}&=&\sqrt{\dfrac{a(1-e^2)}{\mu}}\left[\dfrac{\sin(v+\omega)}{(1+e\cos (v))\sin(i)}\varrho_{N}\right]\vspace{2mm} \\
\dot{\omega}&=&\sqrt{\dfrac{a(1-e^2)}{\mu}}\left[-\dfrac{\cos (v)}{e}\varrho_{R}+\dfrac{(2+e\cos (v))\sin (v)}{e(1+e\cos (v))}\varrho_{T}-\dfrac{\sin(v+\omega) \cot(i)}{1+e\cos (v)}\varrho_{N}\right]\vspace{2mm} \\
\dot{v}&=&\sqrt{\dfrac{\mu}{a^3}}\dfrac{(1+e\cos(v))^2}{(1-e^2)^{3/2}}+\sqrt{\dfrac{a(1-e^2)}{\mu}}\left[\dfrac{\cos (v)}{e}\varrho_{R}-\dfrac{(2+e\cos (v))\sin (v)}{e(1+e\cos (v))}\varrho_{T}\right]
\end{array}
\end{equation}
where $\mu$ is the  gravitational parameter and the overdot symbol denotes the time derivative.

As seen from Eq.~\eqref{gvecoe}, the variational equations for the classical elements are singular at $e=0$ and $i=0$, i.e., for circular and equatorial orbits. The standard method to overcome this issue is to employ a set of modified equinoctial elements defined by
\begin{equation}\label{equinoctialdef}
\begin{array}{l l l}
p&=& {a (1-e^2)}\vspace{1mm}\\
f&=& e \cos(\Omega+\omega)\vspace{1mm}\\
g&=& e \sin(\Omega+\omega)\vspace{1mm}\\
h&=& \tan(i/2) \cos(\Omega)\vspace{1mm}\\
k&=& \tan(i/2) \sin(\Omega)\vspace{1mm}\\
L&=&\Omega+\omega+v
\end{array}
\end{equation}
The dynamics of the modified equinoctial elements are obtained by differentiating \eqref{equinoctialdef} with respect to time and using \eqref{gvecoe}, which results in a highly coupled set of differential equations. In particular, the normal component $\varrho_{N}$ of the perturbing acceleration affects the evolution of all the parameters in \eqref{equinoctialdef}, except for the semiparameter $p$. This is undesirable from the point of view of numerical optimization, which often benefits from sparsity in the system model structure \cite{Junkins19}. In the following, an alternative set of nonsingular elements is presented to address this issue.

\subsection{Ideal Elements}
Let us define the \emph{ideal anomaly} $s$ via the integral relation
\begin{equation}\label{svariable}
s=\int_{t_0}^t\sqrt{\frac{\mu}{a^3(\tau)}}\frac{[1+e(\tau)\cos(v(\tau))]^2}{[1-e^2(\tau)]^{3/2}}\; \text{d}\tau
\end{equation}
where $t$ denotes the actual time and ${t_0}$ is the initial time. Notice from  \eqref{gvecoe} that the time derivative of \eqref{svariable} is equal to $\dot{v}$ for the unperturbed motion ($\varrho_{R}=\varrho_{T}=\varrho_{N}=0$). However, the analogy breaks down in the perturbed case. Based on \eqref{equinoctialdef}-\eqref{svariable}, the following set parameters is proposed in order to describe the orbital motion
\begin{equation}\label{hansenel}
\begin{array}{l l l}
\rho&=& p/R_e \vspace{1mm} \\
e_x&=&f\cos(L-s)+g\sin(L-s) \vspace{1mm} \\
e_y&=&\!\!\!-f\sin(L-s)+g\cos(L-s) \vspace{1mm} \\
h_x&=&h\cos(L-s)+k\sin(L-s) \vspace{1mm} \\
h_y&=&\!\!\!-h\sin(L-s)+k\cos(L-s) \vspace{1mm} \\
\sigma&=&L-s
\end{array}
\end{equation}
where $R_e$ is the radius of the central body (the Earth in this study). The parametrization \eqref{hansenel} is obtained from the modified equinoctial elements by normalizing the semiparameter $p$, subtracting the ideal anomaly from the true longitude to obtain the slowly time-varying parameter $\sigma$, and rotating both the eccentricity vector $[f\;g]^T$ and the ascending node vector $[h\;k]^T$ by an angle $-\sigma$ about the orbit normal. The latter transformation corresponds to rotating the equinoctial frame by an angle $\sigma$ about the orbit normal, so as to match the instantaneous orientation of the \emph{Hansen ideal frame} (see, e.g., \cite{jochim2012significance}). Indeed, $\sigma$ is precisely the in-plane angle between the basis vectors of the Hansel frame and that of the equinoctial frame (in this paper, the initial orientation of the Hansen ideal frame, i.e., the so-called departure frame, is defined by setting $\sigma(0)=L(0)$, so that $s(0)=0$). One can thus refer to the parameters in \eqref{hansenel} as \emph{ideal elements}.

The variational equations for the ideal elements are obtained from \eqref{gvecoe}-\eqref{hansenel} as
\begin{equation}\label{vareq}
\begin{array}{l l l}
\dot{\rho}&=&\dfrac{\sqrt{R_e\rho}}{\sqrt{\mu}\,\kappa_c}\, 2\rho\, \varrho_T \vspace{2mm}\\
\dot{e}_x &=& \dfrac{\sqrt{R_e\rho}}{\sqrt{\mu}\,\kappa_c}\Big\{\;\kappa_c\sin(s)\,\varrho_R+\Big[2\kappa_c \cos(s)+\kappa_s \sin(s)\Big]\varrho_T\Big\} \vspace{2mm}\\
\dot{e}_y &=& \dfrac{\sqrt{R_e\rho}}{\sqrt{\mu}\,\kappa_c}\Big\{-\kappa_c\cos(s)\,\varrho_R+\Big[2\kappa_c\sin(s)-\kappa_s\cos(s)\Big]\varrho_T\Big\} \vspace{2mm}\\
\dot{h}_x &=& \dfrac{\sqrt{R_e\rho}}{\sqrt{\mu}\,\kappa_c}\left[h_x h_y\sin(s)+\kappa_x\cos(s) \right]\varrho_N\vspace{2mm}\\
\dot{h}_y &=& \dfrac{\sqrt{R_e\rho}}{\sqrt{\mu}\,\kappa_c}\left[h_x h_y\cos(s)+\kappa_y\sin(s) \right]\varrho_N\vspace{2mm}\\
\dot{\sigma} &=&\dfrac{\sqrt{R_e\rho}}{\sqrt{\mu}\,\kappa_c}\,\Big[ h_x\sin(s) - h_y\cos(s) \Big]\,\varrho_N
\end{array}
\end{equation}
where
\begin{equation}\label{kappadef}
\begin{array}{l l l}
\kappa_c&=&1+ e_x\cos(s) + e_y\sin(s)\vspace{1mm}\\
\kappa_s&=& e_x\sin(s) - e_y\cos(s)\vspace{1mm}\\
\kappa_x&=& ({1+h_x^2-h_y^2})/{2}\vspace{1mm}\\
\kappa_y&=& ({1-h_x^2+h_y^2})/{2}
\end{array}
\end{equation}
Notice that the right hand side of \eqref{vareq} is independent of the parameter $\sigma$. Moreover, the normal component $\varrho_N$ of the perturbing acceleration affects only the parameters ${h}_x$, ${h}_y$ and $\sigma$.

The particular structure of  \eqref{vareq}-\eqref{kappadef} suggests to adopt the ideal anomaly $s$ as the independent integration variable. Such a procedure is commonly known as regularization \cite{roa2016regularization}. Hereafter, the notation $(\cdot)^\prime=\text{d}(\cdot)/\text{d}s$ will be adopted. From \eqref{equinoctialdef}-\eqref{hansenel} and \eqref{kappadef} one has that
\begin{equation}\label{svar2}
t^\prime=\frac{\text{d}t}{\text{d}s}=\sqrt{\frac{{R_e^3}}{{\mu}}}\,\frac{\rho^{3/2}}{\kappa_c^2}:=\kappa_t
\end{equation}
where the actual time $t$ is now a dependent variable, i.e., $t=t(s)$. Therefore, the regularized dynamics of the ideal elements take on the form
\begin{equation}\label{vareqsd}
[\rho^\prime\;{e}_x^\prime\;{e}_y^\prime\;{h}_x^\prime\;{h}_y^\prime\;\sigma^\prime]^T=\kappa_t\,[\dot{\rho}\;\dot{e}_x\;\dot{e}_y\;\dot{h}_x\;\dot{h}_y\;\dot{\sigma}]^T
\end{equation}
with $[\dot{\rho}\;\dot{e}_x\;\dot{e}_y\;\dot{h}_x\;\dot{h}_y\;\dot{\sigma}]^T$ as in \eqref{vareq}-\eqref{kappadef}.

\section{System Model}\label{sysmod}
In this Section, the proposed parametrization of the orbital motion is exploited to define the system dynamic model used for trajectory optimization. The model includes thrusting, Earth oblateness, and solar eclipse effects.

Let us define the state vector $\mathbf{x}=[\rho\;e_x\;e_y\;h_x\;h_y\;\sigma]^T$ and the acceleration vector $\boldsymbol{\varrho}=[\varrho_R\;\varrho_T\;\varrho_N]^T$. The vector $\boldsymbol{\varrho}$ accounts for the perturbing acceleration $\boldsymbol{\varrho}_{J2}$ induced by the zonal harmonic J2 of the gravitational potential and for a control acceleration $\mathbf{u}$, i.e.,
\begin{equation}
\boldsymbol{\varrho}=\boldsymbol{\varrho}_{J2}+\mathbf{u}
\end{equation}
The contribution $\boldsymbol{\varrho}_{J2}$ can be expressed in terms of $\mathbf{x}$ and $s$  as follows
\begin{equation}\label{j2contrib}
\boldsymbol{\varrho}_{J2}(\mathbf{x},s)=-\frac{3\mu\, \kappa_c^4 J_2 }{2 \, R_e^2\, \rho^4}
\left[
\begin{array}{c}
1 - \dfrac{12(h_x\sin(s)-h_y\cos(s))^2}{(1+ h_x^2 + h_y^2)^2}
\\[4mm]
\dfrac{8(h_x\sin(s)\!-\!h_y\cos(s))(h_x\cos(s)\!+\!h_y\sin(s))}{(1+ h_x^2 + h_y^2)^2}
\\[4mm]
\dfrac{4(h_x\sin(s)-h_y\cos(s))(1-h_x^2-h_y^2)}{(1 +h_x^2 + h_y^2)^2}
\end{array}
\right]
\end{equation}
where $J_2$ is the J2 harmonic coefficient \cite{vallado2001fundamentals}. The regularized dynamics of the state vector $\mathbf{x}$ are obtained from \eqref{vareq}-\eqref{j2contrib} as follows
\begin{equation}\label{pertreldyn}
\mathbf{x}^\prime=\mathbf{f}(\mathbf{x},s)+\mathbf{G}(\mathbf{x},s)\mathbf{u}
\end{equation}
where $\mathbf{f}(\mathbf{x},s)=\mathbf{G}(\mathbf{x},s)\boldsymbol{\varrho}_{J2}(\mathbf{x},s)$ and
\begin{equation*}\setlength\arraycolsep{2pt}
\!\!\mathbf{G}(\mathbf{x},s)\!=\!\!
\dfrac{R_e^2\,\rho^2}{\mu\,\kappa_c^3}\!\!
\left[
\begin{array}{c c c}
0& {2\rho} &0 \vspace{2mm}\\
{\kappa_c\sin(s)}& {2\kappa_c\cos(s)}\!+\!\kappa_s\sin(s)&0\vspace{2mm}\\
\!{-\kappa_c\cos(s)}&{2\kappa_c\sin(s)}\!-\!\kappa_s\cos(s) &0 \vspace{2mm}\\
0& 0 & {h_x h_y\sin(s)\!+\!\kappa_x\cos(s)} \vspace{2mm}\\
0& 0& {h_x h_y\cos(s)\!+\!\kappa_y\sin(s)} \vspace{2mm}\\
0& 0& {h_x\sin(s) - h_y\cos(s)}
\end{array}\right]
\vspace{2mm}
\end{equation*}

The input vector $\mathbf{u}$ in \eqref{pertreldyn} describes the acceleration generated by the spacecraft propulsion system, which in this work consists of a steerable thruster. More specifically, the three-dimensional vector $\mathbf{u}$ is parameterized as
\begin{equation}\label{cinpvector}
\mathbf{u}=\frac{T_{max}}{m}\,\zeta\,\mathbf{q}
\end{equation}
where $T_{max}$ is the maximum deliverable thrust, $m$ is the satellite mass, $\zeta\in[0,\, 1]$ is the throttle factor, and $\mathbf{q}$ is the unitary thrust direction vector ($\| \mathbf{q}\|=1$).  In this setting, the actual control commands are $\zeta$ and $\mathbf{q}$.

The mass variation due to thrusting is accounted for by augmenting system \eqref{pertreldyn}-\eqref{cinpvector} with the regularized mass flow rate equation
\begin{equation}\label{fuel}
m^\prime=-\kappa_t\,\dfrac{m\, \|\mathbf{u}\|}{g_0\, I_{sp}}=
- \kappa_t\,\dfrac{T_{max} }{g_0\, I_{sp}}\,\zeta
\end{equation}
where $g_0$ denotes the standard gravity and $I_{sp}$ is the specific impulse. The dynamic model description is completed by the timing equation \eqref{svar2}.
The full system state vector is then $[\mathbf{x}^T\;m\;t]^T$ and the corresponding regularized dynamic equations are given by \eqref{svar2}, \eqref{pertreldyn} and \eqref{fuel}.

\subsection{Eclipse Conditions}
A spacecraft is shadowed from the Sun when the angle $\theta_{es}$ between the Earth and the Sun, seen from the spacecraft, is smaller than the apparent angular size $\theta_e$ of the Earth radius plus the apparent angular size $\theta_s$ of the Sun radius. A standard definition for the shadow function is then
\begin{equation}\label{shadowf}
\psi=\left\{
\begin{array}{l c c}
0&\text{if}&\theta_{es}\leq \theta_e+\theta_s\\
1&\text{if}&\theta_{es}> \theta_e+\theta_s
\end{array}
\right.
\end{equation}
where $\psi=0$ indicates that that the satellite is shadowed, while no shadowing takes place for $\psi=1$. The angles $\theta_{es}$, $\theta_{e}$ and $\theta_{s}$ are given by
\begin{equation}
\begin{array}{c c l}
\theta_{es}&=&\arccos\left(\dfrac{\mathbf{r}_{es}^T\mathbf{r}_{ss}}{\| \mathbf{r}_{es}\|\, \| \mathbf{r}_{ss}\|}\right)\\[5mm]
\theta_{e}&=&\arcsin(R_e/\|\mathbf{r}_{es}\|)\\[3mm]
\theta_{s}&=&\arcsin(R_s/\|\mathbf{r}_{ss}\|)
\end{array}
\end{equation}
where $\mathbf{r}_{es}$ and $\mathbf{r}_{ss}$ are the Earth and the Sun position vectors relative to the spacecraft, and $R_s$ denotes the Sun radius. The satellite-Sun vector $\mathbf{r}_{ss}$ is given by
\begin{equation}\label{rssdef}
\mathbf{r}_{ss}=\mathbf{r}_{s}(t)+\mathbf{r}_{es}
\end{equation}
where $\mathbf{r}_{s}(t)$ is the ECI position vector of the Sun, which is available as a function of time from ephemeris data. The position vector $\mathbf{r}_{es}$ can be expressed in terms of $\mathbf{x}$ and $s$ as follows
\begin{equation}\label{resdef}
\mathbf{r}_{es}=\kappa_r
\left[
\begin{array}{c}
\cos(s-\sigma)(h_x^2-h_y^2) + 2\sin(s-\sigma) h_x h_y + \cos(s+\sigma)\\[1mm]
\sin(s-\sigma)(h_y^2-h_x^2) + 2\cos(s-\sigma) h_x h_y + \sin(s+\sigma)\\[1mm]
2(h_x\sin(s)-h_y\cos(s))
\end{array}
\right]
\end{equation}
where
\begin{equation*}
\kappa_r=-\frac{R_e \, \rho}{\kappa_c(1+h_x^2+h_y^2)}
\end{equation*}

For the purpose of numerical optimization, the step function \eqref{shadowf} is smoothed by using a logistic function of the form \cite{aziz2019smoothed}
\begin{equation}\label{smoothecl}
\psi_\ell=\frac{1}{1+e^{\,c\, \ell}}
\end{equation}
where $c>0$ is an assignable gain and
$
\ell=\theta_e+\theta_s-\theta_{es}
$.
Notice that such an approximation does not necessarily involve a decrease of modeling accuracy. In fact, eclipse transitions are truly smooth physical phenomena.

\section{Optimal Control Problem}\label{sec4}

In this paper, the control objective is to transfer a satellite from a given initial orbit towards a predefined target orbit in finite time. The target orbit is specified in terms of equinoctial elements (see \eqref{equinoctialdef}) by introducing the reference state vector
\begin{equation}\label{refpar}
\overline{\mathbf{z}}=
\left[
\begin{array}{c}
\overline{p}\\[1mm] \overline{f}\\[1mm] \overline{g}\\[1mm] \overline{h}\\[1mm] \overline{k}
\end{array}
\right]
=
\left[
\begin{array}{c}
\overline{a}(1-\overline{e}^2)\\[1mm] \overline{e} \cos(\overline{\Omega}+\overline{\omega}) \\[1mm] \overline{e} \sin(\overline{\Omega}+\overline{\omega})\\[1mm] \tan(\overline{i}/2)\cos(\overline{\Omega})\\[1mm] \tan(\overline{i}/2)\sin(\overline{\Omega})
\end{array}
\right]
\end{equation}
Let $s_f$ be the final value for the integration variable $s$. According to \eqref{hansenel}, the orbit matching condition at instant $s_f$ can be formalized as
\begin{equation}\label{targetcon}
\mathbf{C}\mathbf{x}(s_f)-\mathbf{D}(\sigma(s_f))\overline{\mathbf{z}}=\mathbf{0}
\end{equation}
where $\mathbf{C}=[\mathbf{I}_{5\times 5}\;\mathbf{0}_{5\times 1}]$ and
\begin{equation}
\mathbf{D}(\sigma)=
\left[
\begin{array}{r r r r r}
1/R_e & 0 & 0 & 0 & 0\\
0 & \cos(\sigma) & \sin(\sigma) & 0 & 0\\
0 & -\sin(\sigma) & \cos(\sigma) & 0 & 0\\
0 & 0 & 0 & \cos(\sigma) & \sin(\sigma)\\
0 & 0 & 0 & -\sin(\sigma) & \cos(\sigma)
\end{array}
\right]
\end{equation}
Thrust limitations due to eclipsing are taken into account in \eqref{cinpvector} by setting
\begin{equation}\label{eclinp}
\zeta=\psi_\ell\, \eta
\end{equation}
and considering $\eta\in[0,\, 1]$ as a new throttle control input.
Substituting  \eqref{eclinp} into \eqref{cinpvector} gives the expression for the control input vector applied to system \eqref{pertreldyn}
\begin{equation}\label{cinpvector2}
\mathbf{u}=\frac{T_{max}}{m}\,\psi_\ell\, \eta\,{\mathbf{q}}
\end{equation}
where the decision variables are $\eta$ and $\mathbf{q}$.

For the system model at hand, the time of flight (TOF), expressed in days, is given by
\begin{equation}
\Delta t=\frac{t(s_f)-t(0)}{86400}
\end{equation}
Moreover, by using \eqref{fuel},\eqref{eclinp}, and the fact that $0\leq \psi_\ell\leq 1$, one can establish the following upper bound on the fuel consumption
\begin{equation}\label{fcub}
\Delta m=\int_{0}^{s_f}\kappa_t\,\frac{T_{max} }{g_0\, I_{sp}}\,\eta\;\,\text{d}s
\end{equation}
where $\kappa_t=\kappa_t(\mathbf{x},s)$ according to \eqref{kappadef}-\eqref{svar2}.
In order to trade off fuel consumption and TOF, we define the performance index
\begin{equation}\label{cfun}
J=(1-\alpha)\,\Delta t +\alpha\, \Delta m
\end{equation}
where $\alpha\in[0,\, 1]$ is a predefined constant. The considered optimal control problem is then formulated as
\begin{equation}\label{optcon}
\begin{aligned}
\underset{\mathbf{q},\,\eta,\,s_f}{\text{min}}\quad &
J=(1-\alpha)\,\Delta t +\alpha\, \Delta m\\[1mm]
\text{s.t.} \quad &
\mathbf{x}(0)=\mathbf{x}_0,\;m(0)=m_0,\;t(0)=t_0  \\[1mm]
& \mathbf{x}^\prime=\mathbf{f}(\mathbf{x},s)+\mathbf{G}(\mathbf{x},s)\,\frac{T_{max}}{m}\,\psi_\ell\,\eta\,\mathbf{q} \\
&m^\prime=\!-\kappa_t\,\frac{T_{max} }{g_0\, I_{sp}}\,\psi_\ell\,\eta\\
&t^\prime\;=\;\kappa_t\\[1mm]
& \mathbf{C}\mathbf{x}(s_f)-\mathbf{D}(\sigma(s_f))\overline{\mathbf{z}}=\mathbf{0}\\[1mm]
&L_{min} \leq\sigma(s_f)+ s_f \leq L_{max}\\[1mm]
&\Delta t_{min} \leq \Delta t \leq \Delta t_{max}\\[1mm]
&0 \leq \eta \leq 1\\[1mm]
&\|\mathbf{q}\|= 1
\end{aligned}
\end{equation}
where $\mathbf{x}_0$, $m_0$ and $t_0$ are given initial conditions and $\{L_{min},\Delta t_{min}\}$, $\{L_{max},\Delta t_{max}\}$ are prescribed lower and upper bounds for the terminal values of the true longitude and of the TOF (one can assign a fixed terminal longitude $\overline{L}$ and a fixed TOF $\Delta \overbar{t\,}$ by setting $L_{min}=L_{max}=\overbar{L}$ and $\Delta t_{min}=\Delta t_{max}=\Delta \overbar{t\,}$). Notice that the final value $s_f$ of the ideal anomaly is itself a decision variable. Problem \eqref{optcon} is in the form of a standard Bolza problem for the nonlinear nonautonomous system defined by the regularized dynamics of the state vector $[\mathbf{x}^T\;m\;t]^T$.

\begin{remark}
The upper bound \eqref{fcub} has been employed to avoid including the stiff function $\psi_\ell=\psi_\ell(\mathbf{x},t,s)$ in the optimization cost. This is done to facilitate numerical optimization. The resulting $\Delta m$ is an excellent approximation of the actual fuel consumption dictated by \eqref{fuel},\eqref{eclinp}, due to the minimization with respect to $\eta$ and the fact that for sufficiently large values of $c$ the shadow function \eqref{smoothecl} closely approximates the step function.
\end{remark}
\begin{remark}
The nonconvex constraint $\|\mathbf{q}\|=1$ has been identified as a major source of computational issues for the solution of problem \eqref{optcon}. For ease of computation, in the software implementation of \eqref{optcon} this constraint is reformulated by introducing the auxiliary decision vector $\mathbf{w}$ and enforcing $\mathbf{q}=\mathbf{w}/\|\mathbf{w} \|$ as well as the convex inequality $\|\mathbf{w}\|^2\leq 1$. The latter is included without loss of generality in order to bound the domain of the OCP.
\end{remark}

A collocation approach is employed to transcribe the continuous-time problem \eqref{optcon} into a static NLP, which is then solved with a suitable nonlinear optimizer.  It is worth recalling that, due to the lack of convexity of problem \eqref{optcon}, the quality of the solution will be highly dependent of the availability of a reasonable initial guess. In the following, a Lyapunov-based guidance strategy addressing the initial guess generation problem is presented.

\section{Initial Guess Generation}\label{initGG}
The initial guess basically provides the initialization point for the NLP solver. In order to generate an initial guess for the solution of problem \eqref{optcon}, one has to determine a suitable value for $s_f$ as well as a candidate trajectory, defined on the interval $s\in[0,\, s_f]$, for the state and control components $\{\mathbf{x}(s),\,m(s),\,t(s),\,\mathbf{q}(s),\,\eta(s)\}$. Technically speaking, the initial guess need not be a feasible solution to problem \eqref{optcon}, therefore any trajectory may be used in principle. However, experience suggests that for large-scale optimization problems such as the one considered herein the initial guess must be reasonably close to feasibility, otherwise the optimizer may fail to return a solution \cite{conway2010}. In this respect, the application of feedback control techniques based on Lyapunov theory has seen a considerable success. A Lyapunov guidance scheme inspired by the Q-law \cite{petropoulos2004low,petropoulos2005refinements} is proposed below. It aims at providing a good initial guess for the optimization of low-thrust orbit transfers involving predefined changes in all orbital elements except for the true anomaly.

\subsection{Lyapunov-based Guidance Scheme}
The first step for the derivation of the guidance scheme is to quantify the deviation between the controlled orbit and the target one (note that the parameters in Section~\ref{sec2} do not quantify such deviation, as they describe the absolute motion). To this aim, we find it convenient to employ a subset of the orbital parameters introduced in \cite{leomanni2020satellite}. In particular, we will make use of the relative inclination $\gamma$ (the angle between the orbital planes of the two orbits), and of the angles $\lambda_2$ and $\lambda_1$ made respectively by the periapsis of the controlled orbit and by the target periapsis, with respect to the relative line of nodes.  By using the results in \cite{leomanni2020satellite} together with \eqref{equinoctialdef}, \eqref{hansenel} and \eqref{refpar}, one can construct a nonlinear coordinate transformation $\mathbf{y}=\mathbf{y}(\mathbf{x},\overline{\mathbf{z}})$ such that
\begin{equation}
\left[
\begin{array}{c c c c c}
a& e& \gamma& \lambda_1& \lambda_2
\end{array}
\right]^T=\mathbf{y}(\mathbf{x},\overline{\mathbf{z}})
\end{equation}
Defining the relative eccentricity vector as
\begin{equation}
\mathbf{e}
=
\left[
\begin{array}{c}
e\cos(\lambda_2)-\overline{e}\cos(\lambda_1)\\
e\sin(\lambda_2)-\overline{e}\sin(\lambda_1)
\end{array}
\right],
\end{equation}
it can be verified that the controlled and the target orbits coincide if $a=\overline{a}$, $\mathbf{e}=\mathbf{0}$ and $\gamma=0$, where $\overline{a}$ and $\overline{e}$ denote respectively the target semimajor axis and eccentricity. As a scalar measure of the deviation of $a$ from $\overbar{a}$, of  $\mathbf{e}$ from $\mathbf{0}$, and of $\gamma$ from $0$, we consider the Lyapunov function candidate
\begin{equation}\label{Lfc}
V(\mathbf{y})=\dfrac{k_a}{2}\sqrt{\dfrac{\mu}{\overline{a}}}\left(\sqrt{\dfrac{\overline{a}}{a}}-1\right)^2+\sqrt{\dfrac{\mu}{ p}}\left[k_e\,\dfrac{\|\mathbf{e}\|^2}{2}+k_\gamma\,\tan^2\!\left(\dfrac{\gamma}{2} \right)
\right]
\end{equation}
where $k_a\geq0$, $k_e\geq0$, $k_\gamma\geq0$ are constant weighting parameters (for notational simplicity, the dependance on constant parameters is not made explicit in the argument of $V$).
The rationale behind the definition of \eqref{Lfc} is similar to that leading to the Q-law. While the Q-law attempts to quantify the \virg{best-case quadratic time-to-go} for the maneuver, the specific form of \eqref{Lfc} attempts to quantify the delta-v needed to bring the actual orbit to the target orbit. Indeed, it can be easily seen that the dimensional unit of \eqref{Lfc} is m/s.

The $s$-derivative of \eqref{Lfc} is evaluated along the trajectory of system \eqref{pertreldyn}, neglecting J2 effects (i.e., $\mathbf{f}(\mathbf{x},s)=\mathbf{0}$). This results in
\begin{equation}\label{dotV}
V^\prime=\frac{\partial V(\mathbf{y})}{\partial\mathbf{y}}\frac{\partial\mathbf{y}(\mathbf{x},\overline{\mathbf{z}})}{\partial \mathbf{x}}\mathbf{G}(\mathbf{x},s)\,\mathbf{u}:=\mathbf{h}^T\,\mathbf{u}
\end{equation}
Substituting \eqref{cinpvector2} into \eqref{dotV} gives
\begin{equation}\label{dotVb}
V^\prime=\mathbf{h}^T\,\frac{T_{max}}{m}\,\psi_\ell\,\eta\,{\mathbf{q}}
\end{equation}
In order to achieve  $V^\prime\leq 0$, the following functional form is assigned to the initial guess of the thrust direction vector $\mathbf{q}$
\begin{equation}\label{qcom}
\mathbf{q}=-\frac{\mathbf{h}}{\|\mathbf{h}\|}
\end{equation}
Substituting \eqref{qcom} into \eqref{dotVb}, one gets
\begin{equation}\label{dotVc}
V^\prime=-\|\mathbf{h}\|\,\frac{T_{max}}{m}\,\psi_\ell\,\eta
\end{equation}
so that $V^\prime\leq 0$, as expected.

In \eqref{dotVc}, the throttle control input $\eta\in[0,\, 1]$ is yet to be specified. The initial guess of $\eta$ is specified according to the coasting mechanism described in \cite{petropoulos2005refinements}, which employs the \emph{time} derivative of $V$ as an indicator of whether to thrust or coast. The time derivative $\dot{V}$ is obtained from \eqref{svar2} and \eqref{dotVc} as
\begin{equation}\label{dotVtime}
\dot{V}=\frac{V^\prime}{\kappa_t}=-\frac{\|\mathbf{h}\|}{\kappa_t}\,\frac{T_{max}}{m}\,\psi_\ell\,\eta
\end{equation}
At each instant $s$, the maximum value $h_{max}$ and minimum value $h_{min}$ of the gain factor ${\|\mathbf{h}\|}/{\kappa_t}=\|\mathbf{h}(\mathbf{x},s)\|/\kappa_t(\mathbf{x},s)$ in \eqref{dotVtime} are predicted over a future time span whose length is equal to one orbital period. This is done by sweeping $\tau$ over the interval $[s,\; s+2\pi]$ in the expression $\|\mathbf{h}(\mathbf{x},\tau)\|/\kappa_t(\mathbf{x},\tau)$, while holding the system state $\mathbf{x}$ constant and equal to the present value $\mathbf{x}(s)$. A minor amendment of the method in \cite{petropoulos2005refinements} is implemented so as to effectively cope with eclipse effects. It boils down to discarding eclipse periods occurring in the interval $[s,\; s+2\pi]$ when computing $h_{max}$ and $h_{min}$ (in other words, $h_{max}$ and $h_{min}$ are evaluated over non-eclipsed orbital arcs). Then, a thrusting efficiency factor is defined as
\begin{equation}
\xi=\dfrac{{\|\mathbf{h}\|}/{\kappa_t}-h_{min}}{h_{max}- h_{min}}
\end{equation}
A coasting phase is enforced whenever the thrusting efficiency is below a predefined threshold $\xi_{cut}\in[0,\,1]$, by setting $\eta=0$ if $\xi< \xi_{cut}$.  Throttle cut-offs due to eclipsing are accounted for by setting $\eta=\psi_\ell$ when $\xi\geq \xi_{cut}$. Summarizing, one has that
\begin{equation}\label{coastinp}
\eta=
\left\{
\begin{array}{lll}
\psi_\ell&\text{if}&\xi \geq \xi_{cut}\\[3mm]
0&\text{if}&\xi< \xi_{cut}
\end{array}
\right.
\end{equation}
For minimum time problems, it is customary to set $\xi_{cut}=0$, so as to fire the thruster whenever it is possible. For $\xi_{cut}>0$ a trade-off is established between TOF and control effort. For additional details about the coasting mechanism, the reader is referred to \cite{petropoulos2005refinements}.

The initial guess for the state and control components $\{\mathbf{x}(s),\,m(s),\,t(s),\,\mathbf{q}(s),\,\eta(s)\}$ is generated by integrating system \eqref{svar2},\eqref{pertreldyn},\eqref{fuel} with the control input \eqref{cinpvector2},\eqref{qcom},\eqref{coastinp}. Due to the way the feedback control policy is defined, the initial guess gets close to a feasible solution to problem \eqref{optcon} as soon as the Lyapunov function \eqref{Lfc} approaches zero. An interesting feature of the proposed guidance scheme is that scaling all the weighting parameters in \eqref{Lfc} by the same positive constant does not change the control policy, and hence it does not affect the state trajectory guess. This allows one to restrict the domain of definition of $\{k_a,\, k_e,\, k_\gamma\}$ to the positive orthant of the unit sphere, by means of the spherical coordinate transformation
\begin{equation}\label{gainpar}
\begin{array}{l l l}
k_a&=&\cos(\varphi_{az})\cos(\varphi_{el})\\
k_e&=&\sin(\varphi_{az})\cos(\varphi_{el})\\
k_\gamma&=&\sin(\varphi_{el}),
\end{array}
\end{equation}
where the free parameters are $\varphi_{az}\in[0,\, \pi/2]$ and $\varphi_{el}\in[0,\, \pi/2]$. In this way, the tuning of the initial guess can be conveniently cast as a two-dimensional search over $\{\varphi_{az},\varphi_{el}\}$ for minimum time problems or a three-dimensional search over $\{\varphi_{az},\varphi_{el},\xi_{cut}\}$ for problems involving fuel optimization. For the case studies presented in Section \ref{numsim}, we adopted a naive search procedure which consists of initializing all the tuning parameters to zero, increasing $\varphi_{az}$ and $\varphi_{el}$ until a reasonably low TOF is achieved, and then increasing $\xi_{cut}$ until the desired TOF/fuel trade-off is reached. By doing so, we were able to generate a suitable trajectory guess in few minutes. The use of the parametrization in \cite{leomanni2020satellite} is instrumental to this purpose, as it allows one to describe the relative motion via a minimal set of variables, which translates into a reduced number of weighting parameters in \eqref{Lfc}. It is also worth remarking that the control policy \eqref{qcom},\eqref{coastinp} is well defined for circular and equatorial orbits, as opposed to the Q-law formulation.

\subsection{Definition of the Initial Mesh}
In order to solve \eqref{optcon} by collocation, one must specify an initial mesh, i.e., a discretization grid for the initial guess. In \emph{hp}-adaptive collocation methods, the mesh is refined during the solution process either by dividing a mesh segment or by increasing the number of collocation points within the segment. Clearly, the size and thus the complexity of the optimization problem will be proportional to the resolution of the initial mesh. Moreover, as with any mesh refinement method, the performance of the \emph{hp} method does depend upon the structure of the initial mesh. For these reasons, the definition of the initial mesh is a critical step of the solution process.
\begin{figure}[!t]
\centering
\includegraphics[width=0.48\textwidth]{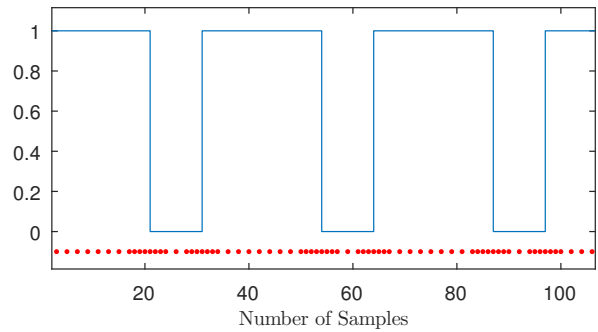}
\caption{Illustration of the initial mesh generation strategy, showing the profile of the shadow function \eqref{shadowf} together with the location the mesh points.}\label{mesh}
\end{figure}
The approach proposed herein for the generation of the initial mesh is as follows. First, system \eqref{svar2},\eqref{pertreldyn},\eqref{fuel},\eqref{cinpvector2},\eqref{qcom},\eqref{coastinp} is numerically integrated and the solution is collected at equally spaced points in the domain $s$, unless an eclipse event is detected. The event detection routine keeps track of all eclipse transitions and associates a sample point to each transition.
The resulting sample sequence is eventually downsampled
by using a method that involves a reduction of the downsample factor in the event proximity. This produces a discretization grid which is finer in correspondence of eclipse transitions (see Fig.~\ref{mesh}), where the system dynamics are inherently stiff. The transitions of the throttle control input $\eta$ arising from \eqref{coastinp} are not treated as events, as they vary less predictably in the solution process and can be assigned freely by the optimizer.

\section{Numerical Simulations}\label{numsim}

In this section, the results of numerical simulations are reported for two benchmark orbit transfer scenarios taken from the literature in order to demonstrate the capabilities of the proposed method. The single-phase OCP \eqref{optcon} is solved using the MATLAB optimal control software GPOPS–II \cite{patterson2014gpops} in combination with the NLP optimizer IPOPT \cite{wachter2006implementation}. GPOPS–II employs an \emph{hp}-adaptive Legendre–Gauss–Radau quadrature orthogonal collocation strategy where the optimal control problem is transcribed into a large sparse NLP, and the NLP is solved on successive mesh iterations until a desired accuracy is achieved. The initial guess for the first mesh iteration is obtained as in Section~\ref{initGG}. In subsequent iterations, the solver is warm-started with the optimal solution from the previous iteration. In this study, the \emph{hp} mesh refinement strategy described in \cite{Liu18hp} is employed. Besides the standard features of \emph{hp} methods, the approach in \cite{Liu18hp} provides the ability to merge mesh segment and to lower the degree of the approximating polynomial, potentially reducing the size of the optimization problem. The mesh refinement accuracy tolerance is set to $\epsilon=5\cdot 10^{-6}$, unless otherwise specified. The number of collocation points per mesh segment is allowed to vary from 4 to 6. The IPOPT optimizer is set up with the linear solver \texttt{ma57} and an error tolerance of $5\cdot 10^{-7}$. The first and second derivatives required by IPOPT are computed analytically by the automatic differentiation software ADiGator \cite{Weinsteinadig}. The constant $c$ in \eqref{smoothecl} is set to $c=298.78$, which according to \cite{aziz2019smoothed} is a realistic value for Earth-centered transfers. All the simulations have been performed on a laptop equipped with a Intel Core I7-5600U processor and 16 GB of RAM.

\subsection{GTO-GEO Transfer}\label{gtogeo}
Herein, we consider the GTO–GEO orbit transfer problem previously solved in \cite{Rao16low,Shannon20lowthrust} via GPOPS-II. The equinoctial elements for the initial GTO and the target GEO are reported in Table \ref{eqgtogeo}.
The GTO elements refer to an orbit with a perigee of 6563.6 km,
an apogee of 42164.3 km, and an inclination of 28.5 deg. The
GEO elements refer to a circular, equatorial orbit with a radius of
42165 km. The Julian date at the beginning of the transfer is $JD_0=2451625.5$.
\begin{table}[!b]
        \caption{Equinoctial elements of the initial and the target orbit for the GTO-GEO transfer}\label{eqgtogeo}
        \centering
        \vspace{3mm}
        \begin{tabular}{c c c c c c}
                \hline
                \hline \\[-2.2ex]
                {Orbit} & p (km)  & f & g& h& k  \\
\hline \\[-2.2ex]
GTO&11359.07&0.7306&0&0.2539676&0\\
GEO&42165&0&0&0&0
\\[0.3mm] \hline
\hline
\end{tabular}
\end{table}
The propulsion parameters are $T_{max}=0.31158$ N and $I_{sp}=1800$ s, while the spacecraft initial mass is $m_0=1200$ kg. All planetary constants for this problem are set equal to those in \cite{Rao16low,Shannon20lowthrust} and are not reported here for brevity.

\begin{figure}[!t]
\centering
\includegraphics[width=0.42\textwidth]{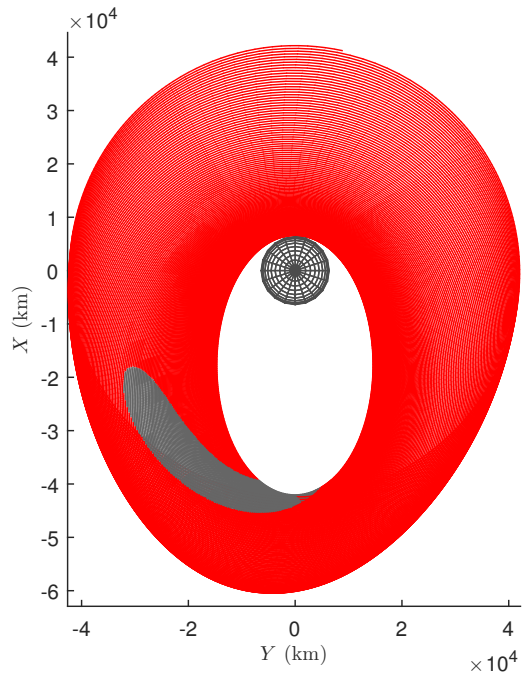}
\caption{Planar projection of the time-optimal GTO-GEO transfer trajectory: thrust phases are colored red, eclipse phases are colored gray.}\label{traj1}
\end{figure}

\begin{figure}[!t]
\centering
\includegraphics[width=0.45\textwidth]{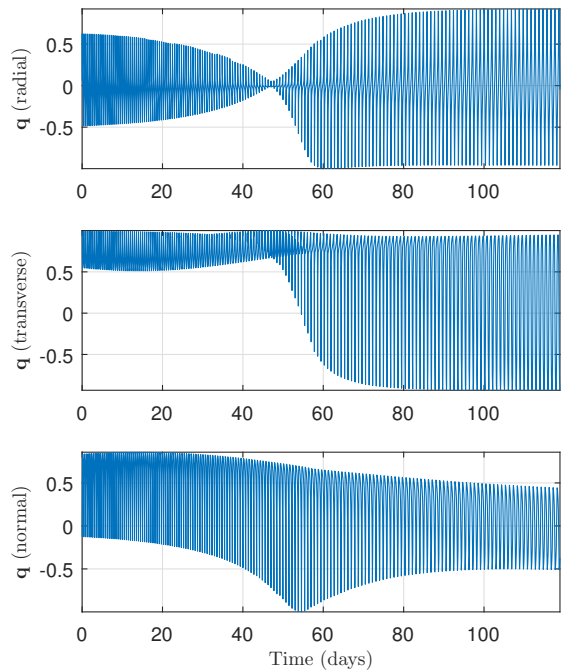}
\caption{Thrust direction vector profile (interpolated from all non-eclipsed trajectory samples) for the time-optimal GTO-GEO transfer.}\label{tvprof}
\end{figure}

Problem \eqref{cfun}-\eqref{optcon} is set up with $\alpha=0$ in order to determine the minimum TOF for the transfer. The initial guess for the problem is generated by tuning the guidance scheme with $\varphi_{az}=7\cdot\pi/180$ rad, $\varphi_{el}=45\cdot\pi/180$ rad and $\xi_{cut}=0$, resulting in trajectory with a TOF of 120.42 days and in a fuel expenditure of 169.79 kg. The trajectory obtained by solving the OCP is depicted in Fig.~\ref{traj1}. It covers approximately 162 revolutions, in terms of the true longitude. The components of the thrust direction vector $\mathbf{q}$ are shown in Fig.~\ref{tvprof}. The throttle control input $\eta$ is always forced to 1 except during eclipse phases. The resulting TOF is 118.74 days and the fuel consumption amounts to 169.38 kg. In this case, they are not far from those provided by the initial guess. The total CPU time for the optimization process is 16 minutes.

The obtained results are compared with those in \cite{Rao16low,Shannon20lowthrust} in Table \ref{mintimesol}. Both the TOF and the fuel consumption are lower than those in~Ref.~\cite{Rao16low} and practically equal to those in Ref.~\cite{Shannon20lowthrust}. Remarkably, the CPU time is one order of magnitude smaller compared to that in Ref.~\cite{Shannon20lowthrust}. This is an especially good figure considering that in \cite{Shannon20lowthrust} the mesh error tolerance is set to $10^{-5}$, while we employed a tolerance value of $5\cdot 10^{-6}$. As a result, we ended up solving a problem with approximately twice the number of NLP variables. Since the main differences between the approach in \cite{Shannon20lowthrust} and the proposed one lie in the parametrization of the dynamic model and in the way the optimization problem is formulated, the observed gain in computational efficiency should be attributed to these two factors.
\begin{table}[!t]
        \caption{Literature comparison for the GTO-GEO transfer}\label{mintimesol}
        \centering
        \vspace{3mm}
        \begin{tabular}{c c c c}
                \hline
                \hline \\[-2.2ex]
                {Source} & TOF  (days) & Fuel (kg) & CPU time (min)  \\
\hline \\[-2.2ex]
Ref.~\cite{Rao16low}&121.22&172.23&n/a\\
Ref.~\cite{Shannon20lowthrust}&118.62&169.44&165 \\
This paper&118.74&169.38&16
\\[0.3mm] \hline
\hline
\end{tabular}
\end{table}
A detailed breakdown of the mesh iteration history is reported in Table \ref{meshiter}. Note that the number of variables in the first iteration is actually higher than that in the last iteration, indicating that the initial mesh is on average more dense than the final one. Although this may seem counter-intuitive, it turned out that initializing the solver with a denser initial mesh does result in a lower number of mesh iterations, which in turns leads to a shorter overall CPU time. Table \ref{meshiter} also demonstrates the ability of the mesh refinement algorithm \cite{Liu18hp} to compress the problem size, thus reducing the computational burden.

\begin{table}[!t]
        \caption{Mesh iteration history for the time-optimal GTO-GEO problem}\label{meshiter}
        \centering
        \vspace{3mm}
        \begin{tabular}{c c c c}
                \hline
                \hline \\[-2.2ex]
                {Iteration} & Mesh error  &NLP var. & CPU time  \\
\hline \\[-2.2ex]
1&9.97e-5&162337&317.7 s\\
2&5.03e-5&127693&95.07 s\\
3&6.15e-6&129901&158.4 s \\
4&5.35e-6&130069&164.1 s\\
5&4.76e-6&130117&183.7 s
\\[0.3mm] \hline
\hline
\end{tabular}
\end{table}

\begin{figure}[!t]
\centering
\includegraphics[width=0.42\textwidth]{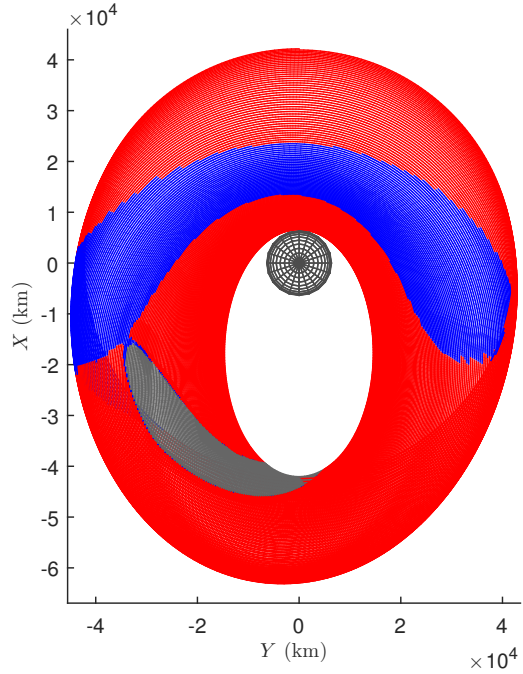}
\caption{Planar projection of the fuel-optimal GTO-GEO transfer trajectory: thrust phases are colored red, eclipse phases are colored gray, and optimal coast phases are colored blue.}\label{traj2}
\end{figure}
\begin{figure}[!t]
\centering
\includegraphics[width=0.45\textwidth]{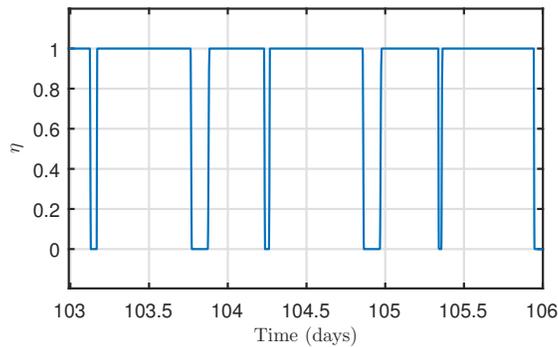}
\caption{Detail of the throttle control input profile for the fuel-optimal GTO-GEO transfer.}\label{thrprof}
\end{figure}

A fuel-optimal GTO-GEO transfer with fixed TOF has been simulated by setting $\alpha=1$ and $\Delta t_{min}=\Delta t_{max}=129.4$ days  in \eqref{cfun}-\eqref{optcon}. The initial guess is left unchanged and equal to that of the minimum time problem, so as to assess the capability of the method to locate optimal coasting arcs without prior indication (i.e., to identify the bang-off-bang structure of the fuel-optimal solution). The problem is solved twice using two different mesh error tolerance levels: $\epsilon=5\cdot 10^{-6}$ and  $\epsilon=5\cdot 10^{-7}$. The solution for $\epsilon=5\cdot 10^{-6}$ features a fuel consumption of 159.39 kg and takes a CPU time of 25.4 min. It involves 7 mesh iterations, with a number of NLP variables in the final mesh equal to 126132. The solution for $\epsilon=5\cdot 10^{-7}$ features a fuel consumption of 159.24 kg and takes a CPU time of 94.6 min. The number of mesh iterations and of NLP variables in the final mesh amount respectively to 11 and 180444. It can be seen that the optimal cost is approximately the same for the two solutions, and that the computational load scales reasonably well with the mesh accuracy. From a qualitative point of view, a lower mesh error tolerance results in a more accurate localization of the coast arcs. The trajectory obtained for $\epsilon=5\cdot 10^{-7}$, reported in Fig.~\ref{traj2}, covers 165 orbital revolutions and cumulates 420 on-off throttle command transitions. A detail of the throttle control input profile is shown in Fig.~\ref{thrprof}.

The mixed time/fuel-optimal problem has been investigated by solving  \eqref{cfun}-\eqref{optcon} with $\alpha=0.5$. In this case, the TOF is not assigned and the optimizer has to search for a pareto-optimal solution. The resulting performance figure is: TOF of 122.6 days, fuel consumption of 162.4 kg and  CPU time of 36.9 min. The pareto-optimal solution displays a TOF that is close to the time-optimal one, while the fuel expenditure is lowered by approximately 7 kg compared to the time-optimal policy. The GTO-GEO transfer results are summarized in Table~\ref{gtogeores}.
\begin{table}[!t]
        \caption{Results for the GTO-GEO transfer}\label{gtogeores}
        \centering
        \vspace{3mm}
        \begin{tabular}{l l c c}
                \hline
                \hline \\[-2.2ex]
                {Solution type} & TOF (days)  & Fuel (kg) & CPU time (min) \\
\hline \\[-2.2ex]
Time-optimal&118.74&169.38&16 \\
Mixed ($\alpha=0.5$)&122.62&162.40 &36.9\\
Fuel-optimal&129.4 (fixed)&159.39&25.4 \\
Fuel-optimal ($\epsilon\!=\!5\cdot 10^{-7}$)&129.4 (fixed)&159.24 &94.6
\\[0.3mm] \hline
\hline
\end{tabular}
\end{table}
\subsection{LEO-GEO Transfer}\label{leogeo}

Herein, we consider the same LEO–GEO orbit transfer scenario studied in \cite{Betts15low,aziz2019smoothed}. The equinoctial elements for the initial LEO and the target GEO are reported in Table \ref{eqleogeo}.
The LEO elements refer to an orbit with an altitude of 500 km above the Earth surface, an inclination of 28.5 deg, and a right ascension of the ascending node of 180 deg. The
GEO elements refer to a circular, equatorial orbit with a radius of about
42241 km. The Julian date at the beginning of the transfer is $JD_0=2457377.5$.
\begin{table}[!b]
        \caption{Equinoctial elements of the initial and the target orbit for the LEO-GEO transfer}\label{eqleogeo}
        \centering
        \vspace{3mm}
        \begin{tabular}{c c c c c c}
                \hline
                \hline \\[-2.2ex]
                {Orbit} & p (km)  & f & g& h& k  \\
\hline \\[-2.2ex]
LEO&6878.140&0&0&-0.2539676&0\\
GEO&42241.095482&0&0&0&0
\\[0.3mm] \hline
\hline
\end{tabular}
\end{table}
The propulsion parameters are $T_{max}=1.445$ N and $I_{sp}=1849.347748$ s, while the spacecraft initial mass is $m_0=1000$ kg. The planetary constants for this problem are set equal to those in \cite{Betts15low,aziz2019smoothed}.

As for the previous case study, we start by evaluating the minimum TOF for the maneuver using $\alpha=0$.  The initial guess for the minimum-time problem is generated by tuning the guidance scheme with $\varphi_{az}=20\cdot\pi/180$ rad, $\varphi_{el}=30\cdot\pi/180$ rad  and $\xi_{cut}=0$. The resulting trajectory converges in approximately 252 orbital revolutions, leading to a time of flight of 44.96 days and a fuel consumption of 289.9 kg. The much higher number of revolutions compared to the GTO-GEO transfer means that we are going to solve a more difficult optimization problem.

\begin{figure*}[!t]
\centering\vspace{2mm}
\includegraphics[width=0.73\textwidth]{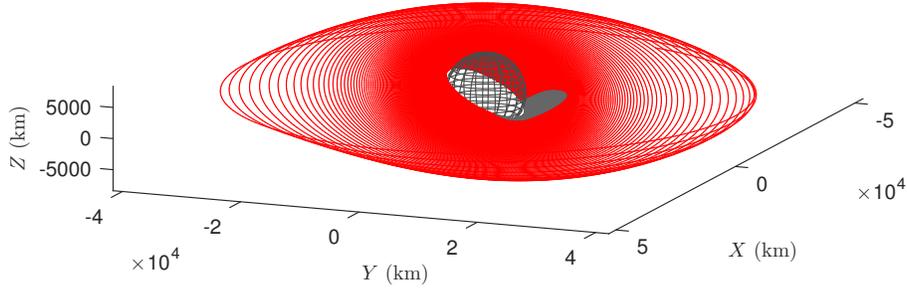}
\caption{Time-optimal LEO-GEO transfer trajectory: thrust phases are colored red, eclipse phases are colored gray.}\label{traj3}
\end{figure*}
\begin{figure}[!t]
\centering
\includegraphics[width=0.45\textwidth]{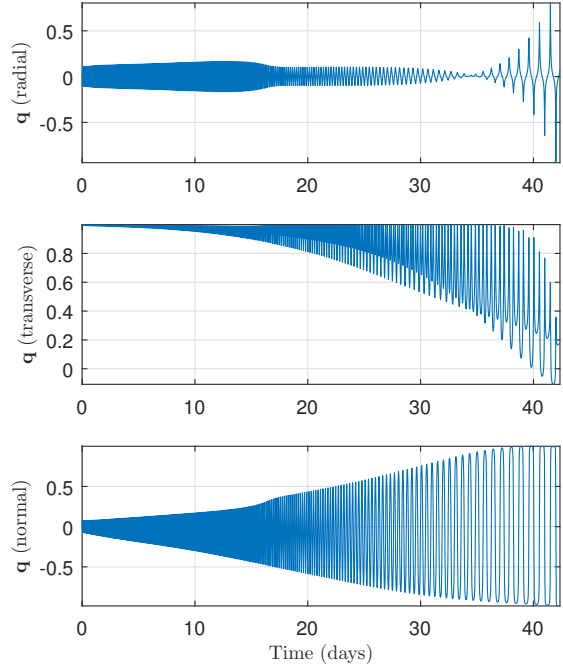}
\caption{Thrust direction vector profile for the time-optimal LEO-GEO transfer.}\label{tvprof2}
\end{figure}
The solution to the time-optimal LEO-GEO problem involves 6 mesh iterations, with a number of NLP variables equal to 252529 in the initial mesh and to 199897 in the final one. The total CPU time taken by the optimization process is 28.3 minutes. The trajectory returned by the optimizer is depicted in Fig.~\ref{traj3}. It covers approximately 251 revolutions in terms of the true longitude. The profile of the thrust direction vector components is shown in Fig.~\ref{tvprof2}. Similarly to the time-optimal GTO-GEO problem, the throttle control input $\eta$ is always forced to 1 except during eclipses. The optimal TOF is 42.37 days and the fuel consumption amounts to 276.70 kg.

\begin{table}[!t]
        \caption{Literature comparison for the LEO-GEO transfer}
        \centering
        \begin{threeparttable}
        \vspace{3mm}
        \begin{tabular}{c c c c}
                \hline
                \hline \\[-2.2ex]
                {Source} & TOF  (days) & Fuel (kg) & CPU time (min)  \\
\hline \\[-2.2ex]
Ref.~\cite{Betts15low}&43.13&281.21&n/a\\
Ref.~\cite{aziz2019smoothed}&44.48&280.88&31.8\tnote{a} \\
This paper&42.37&276.70&28.3\tnote{b}
\\[0.3mm] \hline
\hline
\end{tabular}
\begin{tablenotes}\footnotesize
\item[a] Dual Intel Xeon E5-2860v3, 24 cores
\item[b] Intel Core I7-5600U, single core
\end{tablenotes}
\end{threeparttable}
\label{mintimesol2}
\end{table}

The obtained results are compared with those in \cite{Betts15low,aziz2019smoothed} in Table \ref{mintimesol2}. In these works, the minimum TOF for the maneuver is estimated via heuristic methods and a minimum-fuel problem is solved in which the TOF is fixed and equal to the estimated TOF. As expected, the minimum TOF resulting from the solution to \eqref{cfun}-\eqref{optcon} with $\alpha=0$ is lower than that in \cite{Betts15low,aziz2019smoothed}. Remarkably, the fuel consumption is also lower than that reported in those papers. The total CPU time is close to that obtained in \cite{aziz2019smoothed} using hybrid differential dynamic programming. Note that we employed a J2-perturbed model while Ref.~\cite{aziz2019smoothed} also models the effect of minor orbital perturbations (J3, J4, third-body). However, the results in \cite{aziz2019smoothed} are generated on a dual Intel Xeon E5-2860v3 workstation with a significant part of the workload parallelized over 24 cores, while we employed a laptop processor operated on a single core. Therefore, the computational performance of the proposed solution seems to compare very favorably with that in \cite{aziz2019smoothed}. The CPU time for Ref.~\cite{Betts15low} is not reported, because it lists the total computation time only for the last stage of its multi-phase method (based on the average iteration time reported for the previous stages, we estimated a cumulative CPU time of 47.4 minutes).

\begin{figure*}[!t]
\centering\vspace{2mm}
\includegraphics[width=0.73\textwidth]{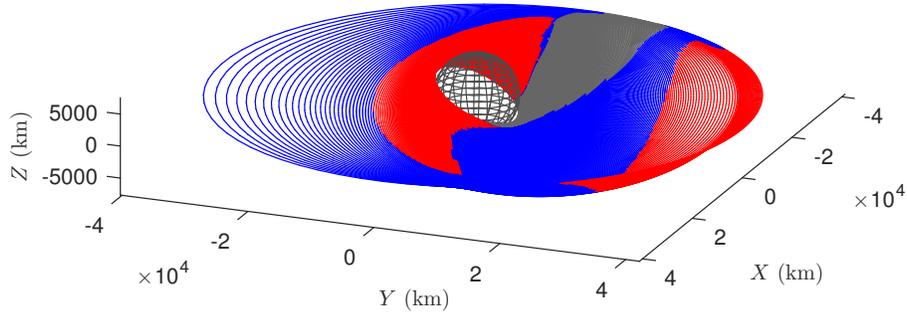}
\caption{Fuel-optimal LEO-GEO transfer trajectory: thrust phases are colored red, eclipse phases are colored gray, and optimal coast phases are colored blue.}\label{traj4}
\end{figure*}
\begin{figure}[!t]
\centering
\includegraphics[width=0.45\textwidth]{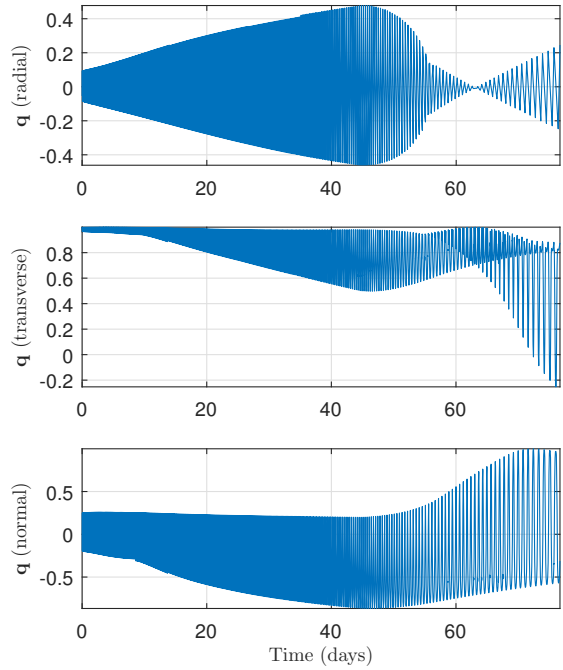}
\caption{Thrust direction vector profile for the fuel-optimal LEO-GEO transfer.}\label{tvprof3}
\end{figure}

Finally, in order to challenge the proposed method, we generated a 405 revolution trajectory guess containing many thrust and coast arcs, by tuning the guidance scheme with $\varphi_{az}=7\cdot\pi/180$ rad, $\varphi_{el}=30\cdot\pi/180$ rad  and $\xi_{cut}=0.3$. This initial guess displays a TOF of $76.7$ days and a fuel consumption of $271.01$ kg. A minimum fuel problem is set up by using $\alpha=1$ and constraining the problem TOF to be equal to the guessed TOF. The problem involves 434304 NLP variables on the initial mesh. Problems of this size are reportedly out of reach even for state-of-the-art low-thrust trajectory optimization software such as Mystic \cite{whiffen2006mystic} (according to \cite{aziz2018low},  computation time limits Mystic to about 250 revolutions for optimized trajectories before switching to the Q-law). We solved the problem twice using the mesh tolerance levels  $\epsilon=5\cdot 10^{-6}$ and  $\epsilon=5\cdot 10^{-7}$. The solution for $\epsilon=5\cdot 10^{-6}$ features a fuel consumption of 236.28 kg and a CPU time of 170 min.
It involves 5 mesh iterations, with a number of NLP variables in the final mesh equal to 290832. The solution for $\epsilon=5\cdot 10^{-7}$ features a fuel consumption of 236.29 kg and a CPU time of 311.5 min. The number of mesh iterations and of NLP variables in the final mesh amount respectively to 9 and 364848. The trajectory obtained for $\epsilon=5\cdot 10^{-7}$, reported in Fig.~\ref{traj4}, covers 405.45 orbital revolutions and cumulates 869 on-off throttle command transitions. The corresponding thrust direction vector profile is shown in Fig.~\ref{tvprof3}.

Note that, for this example, the solution provided by the initial guess is quite far from the optimal one. This demonstrates the capability of the proposed method to effectively explore the solution space. Moreover, it highlights the advantages brought by trajectory optimization with respect to heuristic approaches such as the one in Section~\ref{initGG}, in terms of achievable performance.  The LEO-GEO transfer results are summarized in Table~\ref{leogeores}.

\begin{table}[!t]
        \caption{Results for the LEO-GEO transfer}\label{leogeores}
        \centering
        \vspace{3mm}
        \begin{tabular}{l l c c}
                \hline
                \hline \\[-2.2ex]
                {Solution type} & TOF (days)  & Fuel (kg) & CPU time (min) \\
\hline \\[-2.2ex]
Time-optimal&42.37&276.70&28.3 \\
Fuel-optimal&76.7 (fixed)&236.28&170 \\
Fuel-optimal ($\epsilon\!=\!5\cdot 10^{-7}$)&76.7 (fixed)&236.29 &311.5
\\[0.3mm] \hline
\hline
\end{tabular}
\end{table}

\vspace{-2mm}
\section{Conclusions}\label{sec7}
A direct approach has been presented for the optimization of low-thrust orbit transfer trajectories under eclipse constraints. A specifically conceived parametrization of the orbital motion has been employed in combination with a suitable eclipse smoothing technique in order to define a flexible single-phase optimal control problem formulation.
It has been shown that state-of-the-art pseudospectral collocation algorithms are able to solve this problem effectively. The optimization procedure is complemented by a Lyapunov guidance scheme that can be exploited to generate a reasonable initial guess for the nonlinear solver in short time. The proposed approach is general enough to encompass minimum-time, mimimum-fuel and mixed time/fuel-optimal control problems. Simulations on several relevant missions show that the new parametrization and optimization scheme provide a remarkable improvement in terms of computational efficiency with respect to comparable methods.

\vspace{-2mm}
\bibliographystyle{aiaa-doi}
\bibliography{AIAA20_LTT_minimal}
\end{document}